\documentclass[12pt]{article}
\usepackage{amssymb}
\usepackage{latexsym,bm}
\usepackage{graphicx}
\usepackage{amsmath}
\usepackage{mathrsfs}
\usepackage{epstopdf}

\date{}
\newcounter{mathitem}
\newenvironment{mathitem}
  {\begin{list}{{$(\roman{mathitem})$}}{
   \setcounter{mathitem}{0}
   \usecounter{mathitem}
   \setlength{\topsep}{0pt plus 2pt minus 0pt}
   \setlength{\parskip}{0pt plus 2pt minus 0pt}
   \setlength{\partopsep}{0pt plus 2pt minus 0pt}
   \setlength{\parsep}{0pt plus 2pt minus 0pt}
   \setlength{\leftmargin}{35pt}
   \setlength{\itemsep}{0pt plus 2pt minus 0pt}}}
  {\end{list}}

\begin{document}
\title{The minimum size and maximum diameter of an edge-pancyclic graph of a given order
\footnote{E-mail addresses: {\tt lichengli0130@126.com} (C. Li), {\tt liufeng0609@126.com} (F. Liu), {\tt zhan@math.ecnu.edu.cn} (X. Zhan).}}
\author{\hskip -10mm Chengli Li, Feng Liu and Xingzhi Zhan\thanks{Corresponding author}\\
{\hskip -10mm \small Department of Mathematics, East China Normal University, Shanghai 200241, China}}\maketitle
\begin{abstract}
A $k$-cycle in a graph is a cycle of length $k.$ A graph $G$ of order $n$ is called edge-pancyclic if for every integer $k$ with $3\le k\le n,$ every edge of $G$ lies in a
$k$-cycle. It seems difficult to determine the minimum size $f(n)$ of a simple edge-pancyclic graph of order $n.$ We give lower and upper bounds on $f(n),$ and determine
the maximum diameter of such a graph. In the $3$-connected case, the precise value of $f(n)$ is determined. We also determine the minimum size of a graph of a given order with connectivity conditions in which every edge lies in a triangle.
\end{abstract}


\vskip 8mm

\section{Introduction}

We consider finite simple graphs and use standard terminology and notation from [2] and [9]. The {\it order} of a graph is its number of vertices, and the
{\it size} its number of edges.  A $k$-cycle is a cycle of length $k.$ A graph $G$ of order $n$ is called {\it pancyclic} if for every integer $k$ with $3\le k\le n,$
$G$ contains a $k$-cycle. Putting stronger conditions we have the following two concepts.

{\bf Definition 1.} A graph $G$ of order $n$ is called {\it edge-pancyclic} ({\it vertex-pancyclic}) if for every integer $k$ with $3\le k\le n,$ every edge (vertex)
of $G$ lies in a $k$-cycle.

Clearly, every nonempty edge-pancyclic graph is vertex-pancyclic. Much work has been done on these two classes of graphs ([1], [3], [4], [6], [7], [8], [10]).
Let $g(n)$ denote the minimum size of a vertex-pancyclic graph of order $n.$ Broersma [3] proved that if $n\ge 7$ then $3n/2< g(n)\le 5n/3.$ No conjecture on the
precise value of $g(n)$ is known.

In this paper we give lower and upper bounds on the  minimum size of an edge-pancyclic graph of order $n.$ To do so, we first determine the minimum size of a graph of
order $n$ with connectivity conditions in which every edge lies in a triangle. Finally we determine the maximum diameter of an edge-pancyclic graph of order $n.$

We denote by $V(G)$ and $E(G)$ the vertex set and edge set of a graph $G,$ respectively, and denote by $|G|$ and $e(G)$ the order and size of $G,$ respectively.
Thus $|G|=|V(G)|$ and $e(G)=|E(G)|.$ For a vertex subset $S\subseteq V(G),$ we use $G[S]$ to denote the subgraph of $G$ induced by $S.$ The neighborhood of a vertex
$x$ is denoted by $N(x)$ or $N_G(x)$ and the closed neighborhood of $x$ is $N[x]\triangleq N(x)\cup \{x\}.$ The degree of $x$ is denoted by ${\rm deg}(x).$
We denote by $\delta (G)$ the minimum degree of $G.$
Given two vertex subsets $S$ and $T$ of $G,$ we denote by $[S, T]$ the set of edges having one endpoint in $S$ and the other in $T.$  We denote by $C_n,$ $P_n$ and $K_n$
the cycle of order $n,$ the path of order $n$ and the complete graph of order $n,$ respectively. $\overline{G}$ denotes the complement of a graph $G.$ For two graphs $G$ and $H,$ $G\vee H$ denotes the {\it join} of $G$ and $H,$ which is obtained from the disjoint union $G+H$ by adding edges joining every vertex of $G$ to every vertex of $H.$
An {\it $(x,y)$-path} is a path with endpoints $x$ and $y.$ The distance of two vertices $u$ and $v$ in a graph is denoted by $d(u, v).$

For graphs we will use equality up to isomorphism, so $G=H$ means that $G$ and $H$ are isomorphic.

In Section 2 we determine the minimum size of a $2$-connected or $3$-connected graph of order $n$ in which every edge lies in a triangle, and determine the extremal graphs.
In Section 3 we give lower and upper bounds on the  minimum size of an edge-pancyclic graph of order $n$ and for every integer $k\ge 3$ we construct an
edge-pancyclic graph of order $n=6k^2-5k$ and size $2n-k.$ In Section 4 we determine the maximum diameter of an edge-pancyclic graph of order $n.$

\section{The minimum size of a $2(3)$-connected graph of a given order in which every edge lies in a triangle}

Erd\H{o}s asked the question: What is the minimum size of a connected graph of order $n$ in which every edge lies in a triangle? The answer [5] is $\lfloor (3n-2)/2\rfloor.$
We will re-consider this question with an additional connectivity condition.

We first define four graphs. Let $n\ge 8$ be an integer. If $n$ is even, denote $q=n/2$ and let $C=v_1v_2\ldots v_qv_1$ be a cycle of length $q.$ The graph $A_n$
is obtained from $C$ by adding vertices $u_1,\ldots,u_q$ disjoint from $C$ and adding edges $u_iv_i$ and $u_iv_{i+1}$ for $i=1,\ldots,q$ where $v_{q+1}=v_1.$
Now suppose that $n$ is odd. The next three graphs will be constructed from $A_{n-1}$ whose vertex set is $\{v_1,\ldots, v_{(n-1)/2}\}\cup \{u_1,\ldots,u_{(n-1)/2}\}$ as above.
The graph $F_n$ is obtained from $A_{n-1}$ by adding a vertex $x$ and adding edges $xv_1$ and $xv_2.$ The graph $G_n$ is obtained from $A_{n-1}$ by adding a vertex $y$
and adding edges $yu_1$ and $yv_1.$ The graph $H_n$ is obtained from $A_{n-1}$ by first subdividing the edge $v_1v_2$ using a vertex $z$ and then adding
the edge $zu_1.$ The graphs $A_{12},$ $F_{13},$ $G_{13}$ and $H_{13}$ are depicted in Figure 1.
\begin{figure}[h]
\centering
\includegraphics[width=0.8\textwidth]{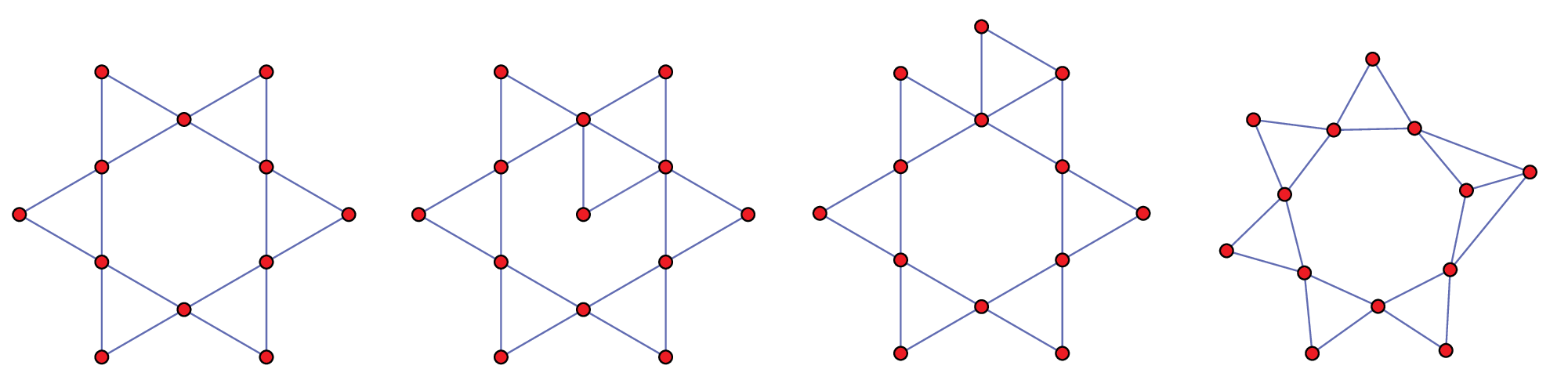}
\caption{The graphs $A_{12}$, $F_{13}$, $G_{13}$ and $H_{13}$}
\end{figure}

{\bf Lemma 1.} {\it Let $G$ be a $2$-connected graph of order $n\ge 6.$ If every edge of $G$ lies in a triangle, then
$$
e(G)\ge \left\lceil\frac{3n}{2}\right\rceil. \eqno (1)
$$
If $n\ge 8$ and $n$ is even, equality holds in (1) if and only if $G=A_n;$ if $n\ge 9$ and $n$ is odd, equality holds in (1) if and only if $G\in \{F_n, G_n, H_n\}.$
}

{\bf Proof.} We use induction on the order $n.$ It is not difficult to verify that Lemma 1 holds for the initial orders $6,$ $7,$ $8,$ $9.$
Next let $n\ge 10$ and we suppose that Lemma 1 holds for graphs of any order $p$ with $6\le p\le n-1.$ We distinguish two cases according as $n$ is even or odd.
Note that $\lceil 3n/2 \rceil=3n/2$ if $n$ is even and $\lceil 3n/2 \rceil=(3n+1)/2$ if $n$ is odd.

Case 1. $n$ is even.

We first prove inequality (1). Since  $G$ is $2$-connected, $\delta(G)\ge 2.$ If $\delta(G)\ge 3,$ then $e(G)\ge n\delta(G)/2\ge 3n/2.$ Suppose that $\delta(G)=2$ and choose
a vertex $v$ with ${\rm deg}(v)=2.$ Let $N(v)=\{x,y\}.$ Since $vx$ lies in a triangle, $x$ and $y$ are adjacent. We assert that $G-v$ is $2$-connected. Otherwise, any cut-vertex
of $G-v$ is also a cut-vertex of $G,$ a contradiction. If the edge $xy$ lies in a triangle in $G-v,$ then by the induction hypothesis we have $e(G-v)\ge (3(n-1)+1)/2=(3n-2)/2.$
Thus $e(G)\ge (3n-2)/2+2>3n/2.$

Next suppose that $xy$ does not lie in any triangle in $G-v;$ i.e., $N(x)\cap N(y)=\{v\}.$ Recall that to {\it contract} an edge $f=ab$ is to replace $a$ and $b$ with a single vertex whose incident edges are the edges other than $f$ that were incident with $a$ or $b.$ Let $H$ be the graph obtained from $G-v$ by contracting the edge $xy$ to a new vertex $z.$
Then $H$ is a connected simple graph in which every edge lies in a triangle. If $H$ is $2$-connected, by the induction hypothesis we have $e(H)\ge 3(n-2)/2,$ implying that
$$
e(G)\ge e(H)+3\ge 3(n-2)/2+3=3n/2.
$$

Now suppose that $H$ is not $2$-connected. Then $z$ is the only cut-vertex of $H,$ since otherwise any other cut-vertex would be a cut-vertex of $G,$ contradicting  the condition that
$G$ is $2$-connected. Thus, the block-cut-vertex tree of $H$ is a star. Let $H_1,\ldots, H_s$ be the blocks of $H$ with $s\ge 2,$ each of which contains $z.$

{\bf Claim.} $|H_i|\ge 5$ and if $|H_i|=5$ then $e(H_i)\ge 8=(3|H_i|+1)/2$ for each $i=1,\ldots,s.$

Choose an arbitrary but fixed index $i$ with $i=1,\ldots,s.$ Denote $H_i^{\prime}=H_i-z.$ Then  $H_i^{\prime}$ is a component of $G-\{x,y\}.$ Since $G$ is $2$-connected, $\{x,y\}$ is a minimum vertex cut of $G.$ Hence $N_{H_i^{\prime}}(x)\neq\emptyset$ and $N_{H_i^{\prime}}(y)\neq\emptyset.$ Let $x_1\in N_{H_i^{\prime}}(x)$ and $y_1\in N_{H_i^{\prime}}(y).$
Since each of $xx_1$ and $yy_1$ lies in a triangle, there exist vertices $x_2,y_2\in V(H_i^{\prime})$ such that $xx_1x_2x$ and $yy_1y_2y$ are triangles. Since $N(x)\cap N(y)=\{v\},$
the four vertices $x_1,x_2,y_1,y_2$ are pair-wise distinct. With the additional vertex $z,$ we deduce that $|H_i|\ge 5.$

Suppose $|H_i|=5.$ Then $[\{x_1, x_2\}, \{y_1, y_2\}]\neq\emptyset,$ since $H_i$ is $2$-connected. Without loss of generality suppose $x_1y_1\in [\{x_1, x_2\}, \{y_1, y_2\}].$
In $G,$ the edge $x_1y_1$ lies in a triangle $px_1y_1p.$ Clearly $p\neq v.$ Since $N(x)\cap N(y)=\{v\},$ $p\notin \{x,y\}.$ Thus $p\in \{x_2, y_2\}.$ It follows that
in addition to the seven edges $zx_1, zx_2, x_1x_2, zy_1, zy_2, y_1y_2, x_1y_1,$ there exists at least one more edge in $H_i.$ Hence $e(H_i)\ge 8.$ This proves the claim.

Observe that $|H|=n-2$ and $\sum_{i=1}^{s}|H_i|=n+s-3.$ We may suppose that the first $r$ blocks $H_1,\ldots,H_r$ have even orders and the remaining blocks $H_{r+1},\ldots,H_s$ have
odd orders. Applying the above claim (for those $j$ with $|H_j|=5$) and the induction hypothesis (for those $j$ with $|H_j|\ge 6$) we have $e(H_i)\ge 3|H_i|/2$ for $1\le i\le r$
and $e(H_i)\ge (3|H_i|+1)/2$ for $r+1\le i\le s.$ Thus
\begin{align*}
e(G)=3+\sum_{i=1}^{s}e(H_i)&\ge 3+\sum_{i=1}^{s}\frac{3|H_i|}{2}+\frac{s-r}{2}\\
                           &=3+\frac{3(n+s-3)}{2}+\frac{s-r}{2}\\
                           &=\frac{3n+3(s-1)+(s-r)}{2}\\
                           &\ge\frac{3n+3}{2}\\
                           &>\frac{3n}{2}.
\end{align*}

Now we determine the extremal graphs. Suppose that $G$ is a $2$-connected graph of order $n$ and size $3n/2$ and each edge of $G$ lies in a triangle. We will prove that
$G=A_n.$

If $\delta(G)\ge 3,$ then $G$ is cubic. Let $u$ be a vertex and let $N(u)=\{x_1, x_2, x_3\}.$ Since every edge lies in a triangle, $e(G[\{x_1,x_2,x_3\}])\ge 2.$ Without loss of generality, assume that $x_1$ and $x_2$ are adjacent and $x_2$ and $x_3$ are adjacent. If $x_1$ and $x_3$ are adjacent, then $V(G)=\{u,x_1,x_2,x_3\},$ contradicting $n\ge 10.$ Let $x_4$ be the neighbor of $x_1$ other than $u, x_2.$ Then the edge $x_1x_4$ lies in a triangle $x_1x_4x_5x_1.$ Since $x_4,x_5\notin \{u,x_2\},$ we have ${\rm deg}(x_1)\ge 4,$ a contradiction. Hence $\delta(G)=2.$

Let $v$ be a vertex of degree $2$ and let $N(v)=\{x,y\}.$ Combining the condition $e(G)=3n/2$ and the above proof of (1) we deduce that
(i) $x$ and $y$ are adjacent; (ii) $G-v$ is $2$-connected; (iii) $N(x)\cap N(y)=\{v\}.$  Let the graph $H$ and the vertex $z$ be defined as above. 
Clearly $H$ is a graph of order $n-2$ and size $3(n-2)/2$ in which every edge lies in a triangle. We have proved above that if $H$ is not $2$-connected then
$e(G)>3n/2,$ contradicting our assumption that $e(G)=3n/2.$ Hence $H$ is $2$-connected. By the induction hypothesis, $H=A_{n-2}.$ The degree of a vertex in $H$ 
is either $2$ or $4.$ If ${\rm deg}(z)=2,$ then either $G$ is not $2$-connected or $G$ contains three edges which do not lie in any triangle, a contradiction. 
Thus ${\rm deg}(z)=4.$

Denote $m=(n-2)/2$ and let $V(H)=\{v_1,\ldots,v_{m},u_1,\ldots,u_m\}$ where ${\rm deg}(v_i)=4$ and ${\rm deg}(u_i)=2$ for $1\le i\le m.$ Without loss of generality, suppose
that $z=v_2.$ Thus $N_H(z)=\{v_1,v_3,u_1,u_2\}.$ In $G,$ every vertex in $\{v_1,v_3,u_1,u_2\}$ has exactly one neighbor in $\{x, y\}.$ Further, using the conditions
that in $G,$ each of $v_1u_1$ and $v_3u_2$ lies in a triangle and $G$ is $2$-connected, we deduce that $G=A_n.$

Conversely, $A_n$ is a $2$-connected graph of order $n$ and size $3n/2$ in which every edge lies in a triangle.

Case 2. $n$ is odd.

In this case the proof ideas are almost the same as in Case 1, but there are subtle differences in counting. For completeness we give the details. We first prove inequality (1); i.e.,
$$
e(G)\ge \frac{3n+1}{2}. \eqno (2)
$$
Since  $G$ is $2$-connected, $\delta(G)\ge 2.$ If $\delta(G)\ge 3,$ then $e(G)\ge (3(n-1)+4)/2=(3n+1)/2,$ where we have used the fact that $G$ cannot be cubic since $n$ is odd.

Suppose that $\delta(G)=2$ and choose a vertex $v$ with ${\rm deg}(v)=2.$ Let $N(v)=\{x,y\}.$ Since $vx$ lies in a triangle, $x$ and $y$ are adjacent. We assert that $G-v$ is $2$-connected, since $G$ is so. If the edge $xy$ lies in a triangle in $G-v,$ then by the induction hypothesis we have $e(G-v)\ge 3(n-1)/2.$
Thus $e(G)\ge 3(n-1)/2+2=(3n+1)/2.$

Next suppose that $xy$ does not lie in any triangle in $G-v;$ i.e., $N(x)\cap N(y)=\{v\}.$ Let $H$ be the graph obtained from $G-v$ by contracting the edge $xy$ to a new vertex $z.$
Then $H$ is a connected simple graph in which every edge lies in a triangle. If $H$ is $2$-connected, by the induction hypothesis we have $e(H)\ge (3(n-2)+1)/2,$ implying that
$$
e(G)\ge e(H)+3\ge (3(n-2)+1)/2+3=(3n+1)/2.
$$

Now suppose that $H$ is not $2$-connected. Then $z$ is the only cut-vertex of $H.$ Thus, the block-cut-vertex tree of $H$ is a star. Let $H_1,\ldots, H_s$ be the blocks of $H$ with $s\ge 2,$ each of which contains $z.$ The claim in Case 1 still holds.

Observe that $|H|=n-2$ and $\sum_{i=1}^{s}|H_i|=n+s-3.$ We may suppose that the first $r$ blocks $H_1,\ldots,H_r$ have even orders and the remaining blocks $H_{r+1},\ldots,H_s$ have
odd orders. Applying the above claim (for those $j$ with $|H_j|=5$) and the induction hypothesis (for those $j$ with $|H_j|\ge 6$) we have $e(H_i)\ge 3|H_i|/2$ for $1\le i\le r$
and $e(H_i)\ge (3|H_i|+1)/2$ for $r+1\le i\le s.$ Thus, as in Case 1,
$$
e(G)=3+\sum_{i=1}^{s}e(H_i)\ge \frac{3n+3}{2}>\frac{3n+1}{2}.
$$
In any case, inequality (2) is proved.

Now we determine the extremal graphs. Suppose that $G$ is a $2$-connected graph of order $n$ and size $(3n+1)/2$ and each edge of $G$ lies in a triangle. We will prove that
$G\in \{F_n, G_n, H_n\}.$

We assert that  $\delta(G)=2.$ If $\delta(G)\ge 3,$ then the degree sequence of $G$ is $4,3,3,\ldots,3.$ Let $u$ be a vertex of degree $3$ and let $N(u)=\{x_1, x_2, x_3\}.$ Since every edge lies in a triangle, $e(G[\{x_1,x_2,x_3\}])\ge 2.$ Without loss of generality, assume that $x_1$ and $x_2$ are adjacent and $x_2$ and $x_3$ are adjacent. Recall that we have assumed $n\ge 10.$
We assert that $x_1$ and $x_3$ are nonadjacent. Otherwise $N[u]$ is a $4$-clique and one of $x_1,x_2,x_3$ must have degree $4$ and it is a cut-vertex, contradicting the condition
that $G$ is $2$-connected.

Let $x_4$ be a neighbor of $x_1$ other than $u, x_2.$ Since the edge $x_1x_4$ lies in a triangle, either ${\rm deg}(x_1)=4$ or $x_4$ and $x_2$ are adjacent.
Similarly, let $x_5$ be a neighbor of $x_3$ other than $u, x_2.$ Then either ${\rm deg}(x_3)=4$ or $x_5$ and $x_2$ are adjacent. Thus, among the three vertices
$x_1,$ $x_2,$ and $x_3,$ either one has degree $5,$ or two have degree $4,$ which is a contradiction in either case. Hence $\delta(G)=2.$

Let $v$ be a vertex of degree $2$ and let $N(v)=\{x,y\}.$ Then (i) $x$ and $y$ are adjacent; (ii) $G-v$ is $2$-connected.

First suppose that the edge $xy$ lies in a triangle of $G-v.$  We have
$$
e(G-v)=e(G)-2=\frac{3n+1}{2}-2=\frac{3(n-1)}{2}.
$$
By the induction hypothesis, $G-v=A_{n-1},$ implying that $G\in\{F_n, G_n\}.$

Now suppose that the edge $xy$ does not lie in any triangle of $G-v.$ Let $R$ be the graph obtained by contracting the edge $xy$ to a new vertex $z.$ Then $R$ is a
$2$-connected graph of order $n-2$ and size $(3(n-2)+1)/2$ in which every edge lies in a triangle. By the induction hypothesis, $R\in\{F_{n-2},G_{n-2},H_{n-2}\}.$
Using the three properties of $G:$ (i) $x$ and $y$ have no common neighbor in $G-v;$ (ii) every edge of $G$ lies in a triangle; (iii) $G$ is $2$-connected, we deduce
the following conclusions. If $R=F_{n-2}$ then ${\rm deg}_R(z)\in\{4, 5\}$ and $G=F_n;$ if $R=G_{n-2}$ then ${\rm deg}_R(z)\in\{4, 5\}$ and $G=G_n;$ if $R=H_{n-2}$
then ${\rm deg}_R(z)=4$ and $G=H_n.$

Conversely, each of $F_n, G_n,$ and $H_n$ is a $2$-connected graph of order $n$ and size $(3n+1)/2$ in which every edge lies in a triangle. This completes the proof. \hfill$\Box$

Contracting an edge in a simple graph can produce multiple edges. By deleting from a general graph $R$ all loops and, for every pair of adjacent vertices,
all but one edge joining them, we obtain the {\it underlying simple graph} of $R.$

We denote by $W_n$ the wheel graph of order $n;$ i.e., $W_n=K_1\vee C_{n-1}$ where $\vee$ means join.

{\bf Lemma 2.} {\it Let $G$ be a $3$-connected graph of order $n$ in which every edge lies in a triangle. Then
$$
e(G)\ge 2n-2. \eqno (3)
$$
Equality holds if and only if $G=W_n.$  }

{\bf Proof.} We use induction on the order $n.$ Since $G$ is $3$-connected, $n\ge 4$ and $\delta(G)\ge 3.$ It is easy to verify that the cases $n=4,\, 5$ hold.
In fact, $K_4$ is the only $3$-connected graph of order $4,$ and there are exactly three $3$-connected graphs of order $5$ in which every edge lies in a triangle.

Now we let $n\ge 6$ and suppose that Lemma 2 holds for all graphs of orders less than $n.$ If $\delta(G)\ge 4,$ then $e(G)\ge \delta(G)n/2\ge 4n/2=2n>2n-2,$ showing (3).

Next suppose $\delta(G)=3.$ Choose a vertex $v$ of degree $3$ and let $N(v)=\{x, y, z\}.$ Since every edge lies in a triangle, $e(G[\{x, y, z\}])\ge 2.$ Without loss of generality,
suppose that $xy,\, yz\in E(G)$. Let $R$ be the graph obtained from $G$ by contracting the edge $vx$ to a vertex  $v^{\prime}$ and let $H$ be the underlying simple graph
of $R.$ Then $e(H)=e(G)-2$ if $x$ and $z$ are nonadjacent and $e(H)=e(G)-3$ otherwise. Clearly every edge of $H$ lies in a triangle.

{\bf Claim.} $H$ is $3$-connected.

To the contrary, suppose that $H$ is not $3$-connected. Then $H$ contains a vertex cut $S$ of cardinality $2.$ Since $G$ is $3$-connected, $S$ is not a vertex cut of $G$
and hence $v^{\prime}\in S.$ Let $S=\{v^{\prime}, w\}.$ Then $\{v, x, w\}$ is a minimum vertex cut of $G.$ Let $G_1$ and $G_2$ be two distinct components of $G-\{v, x, w\}.$
We have $N_{G_1}(v)\neq\emptyset$ and $N_{G_2}(v)\neq\emptyset.$ Thus either $y\in V(G_1)$ and $z\in V(G_2),$ or $y\in V(G_2)$ and $z\in V(G_1),$ contradicting
the fact that $y$ and $z$ are adjacent. This proves the claim.

By the induction hypothesis, $e(H)\ge 2(n-1)-2=2n-4.$ Hence $e(G)\ge e(H)+2\ge 2n-2.$ This proves (3).

Clearly the wheel graph $W_n$ is an extremal graph for the minimum size $2n-2.$ Next we show that it is the unique extremal graph. Let $G$ be a $3$-connected graph of order $n$
and size $2n-2$ in which every edge lies in a triangle. We have $\delta(G)=3.$ Choose a vertex $v$ of degree $3$ and let $N(v)=\{x, y, z\}.$ Define the vertex $v^{\prime}$
and the graph $H$ as above. We have $e(H)\ge 2n-4.$ Then the condition $e(G)=2n-2$ implies that $e(H)=2n-4$ and the two vertices $x$ and $z$ are nonadjacent.

By the induction hypothesis, $H=W_{n-1}.$ The condition $n-1\ge 5$ implies that the central vertex of $H$ has degree at least $4$ and all other vertices have degree $3.$
We assert that $v^{\prime}$ is not the central vertex. Otherwise let $q$ be the vertex such that $y, z, q$ appear consecutively on the $(n-2)$-cycle of $H.$ Since
$N(v)=\{x, y, z\}$ and the two vertices $x$ and $z$ are nonadjacent, $x$ and $q$ are adjacent and the $4$-cycle $xvzqx$ is chordless. Consequently the edge $zq$ does not
lie in any triangle in $G,$ a contradiction. Thus $v^{\prime}$ lies on the $(n-2)$-cycle of $H.$

Let $N_H(v^{\prime})=\{c,f,g\}$ where $c$ is the central vertex of $H.$ Then $\{y,z\}\subset \{c,f,g\}.$ Since $y$ and $z$ are adjacent and $f$ and $g$ are nonadjacnet,
$\{f,g\}\neq \{y,z\}.$ Hence $c\in \{y,z\}.$ We assert that $c\neq z.$ Otherwise if $f=y$ then the edge $xg$ does not lie in any triangle in $G,$ and if $g=y$ then the
edge $xf$ does not lie in any triangle in $G,$ where we have used the fact that $x$ and $z$ are nonadjacent. Thus $c=y.$ It follows that either $f=z$ or $g=z.$  In both cases
it is easy to see that $G=W_n.$ This completes the proof. \hfill$\Box$

\section{Bounds on the minimum size of an edge-pancyclic graph of a given order}

An immediate consequence of Lemma 2 is the following result.

{\bf Theorem 3.} {\it Let $G$ be a $3$-connected edge-pancyclic graph of order $n.$ Then $e(G)\ge 2n-2$ and equality holds if and only if $G=W_n.$ }

Combining Lemma 1 and the fact that the wheel $W_n$ is an edge-pancyclic graph of order $n$ and size $2n-2,$ we obtain the following result.

{\bf Theorem 4.} {\it Let $f(n)$ be the minimum size of an edge-pancyclic graph of order $n\ge 4.$ Then
$$
\frac{3n}{2}\le f(n) \le 2n-2.
$$
}

We remark that the cases $n=4, 5$ of Theorem 4 can be verified directly. $K_4$ is the only edge-pancyclic graph of order $4,$ and there are three edge-pancyclic graphs of order
$5$ with sizes $8,$ $9,$ and $10,$ respectively.

A computer search yields that $f(n)=2n-2$ for $4\le n\le 12.$ One might conjecture that this is the case for general orders $n.$ Our next result shows that this is far from
the truth.

{\bf Theorem 5.} {\it Given any integer $k\ge 3,$ let $n=6k^2-5k.$ Then there exists an edge-pancyclic graph of order $n$ and size $2n-k.$ }

{\bf Proof.} The {\it fan graph} of order $s,$ denoted by $F_s,$ is the join of a vertex $K_1$ and a path $P_{s-1}.$ Thus $F_s=K_1\vee P_{s-1}.$ Let $A$ and $B$ be
two vertex-disjoint graphs both isomorphic to $F_{3k-2}$ with $V(A)=\{v,v_1,v_2,\ldots, v_{3k-3}\}$ and $V(B)=\{u,u_1,u_2,\ldots, u_{3k-3}\}$ where $v$ and $u$ are the central vertices of $A$ and $B,$ respectively. Let $H(k)$ be the graph obtained from $A$ and $B$ by adding edges $vu,$ $vu_{3k-3}$ and $uv_{3k-3}.$ See Figure 2 for an illustration
of $H(k).$ Note that $H(k)$ has order $6k-4$ and size $12k-11.$
\begin{figure}[h]
\centering
\includegraphics[width=0.6\textwidth]{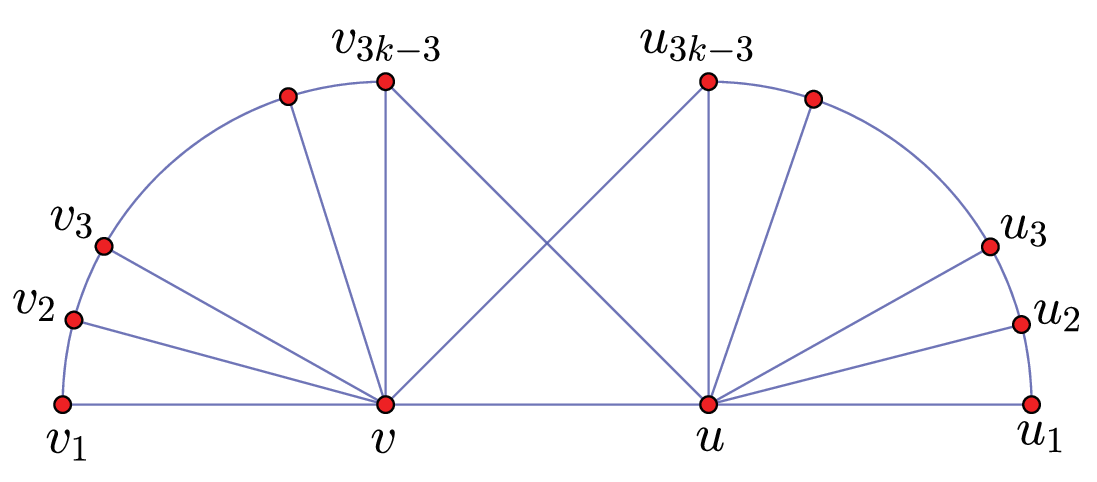}
\caption{The graph $H(k)$}
\end{figure}

Now we construct a graph $G(k)$ as follows. Take a $k$-cycle $C_k=x_1x_2\cdots x_kx_1.$ Make the convention $x_{k+1}=x_1.$ For every $i$ with $1\le i\le k,$ replace the edge
$x_ix_{i+1}$ by a copy of $H(k)$ to obtain the graph $G(k)$ where we identify the vertices $x_i$ and $x_{i+1}$ by $v_1$ and $u_1,$ respectively. The graph $G(3)$ has order
$39$ and it is depicted in Figure 3.
\begin{figure}[h]
\centering
\includegraphics[width=0.45\textwidth]{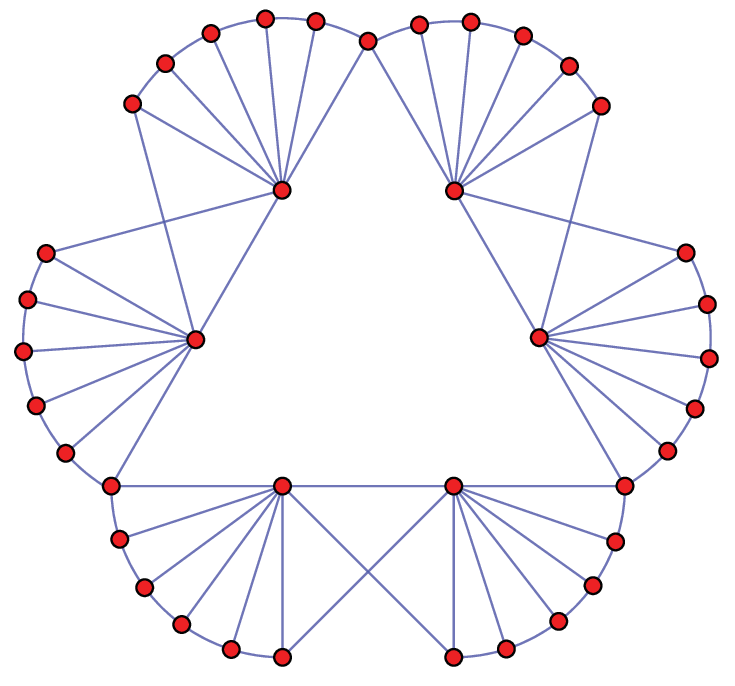}
\caption{The graph $G(3)$}
\end{figure}

Clearly $G(k)$ has order $n=6k^2-5k$ and size $2n-k.$ Next we verify that it is edge-pancyclic.

For integers $a$ and $b$ with $a\le b,$ we use $[a, \,b]$ to denote the set of those integers $c$ satisfying $a\le c\le b.$
An edge $f$ of a graph $R$ of order $r$ is said to be {\it pancyclic} if for every integer $q\in [3,\, r],$ $f$ lies in a $q$-cycle of $R.$

Denote $\Omega=\{vv_i|\, 1\le i\le 3k-3\}\cup \{v_iv_{i+1}|\, 1\le i\le 3k-4\}\cup \{vu, vu_{3k-3}\},$ a subset of the edge set of $H(k).$ The graph $H(k)$ has the following properties:
\begin{mathitem}
\item For every integer $p\in [3,\, 6k-5],$ $H(k)$ contains a $(v_1,u_1)$-path of length $p.$

\item Every edge in $\Omega$ lies in  a $(v_1,u_1)$-path of length $6k-5.$

\item For every index $i\in [1,\, 3k-4],$ each edge $v_iv_{i+1}$ is pancyclic in $H(k),$ and $v_iv_{i+1}$ lies in a $(v_1,u_1)$-path of length $3k-i+1.$

\item $vu_{3k-3}$ is pancyclic in $H(k),$ and $vu_{3k-3}$ lies in a $(v_1,u_1)$-path of length $4.$

\item For every integer $p\in [3,\, 3k-1],$ $vu$ lies in $p$-cycle, and $vu$ lies in a $(v_1,u_1)$-path of length $3.$

\item For every index $i\in [1,\, 3k-3]$ and for every integer $p\in [3,\, 6k-i-3],$ each edge $vv_i$ lies in a $p$-cycle, and $vv_i$ lies in
 a $(v_1,u_1)$-path of length $3k-i+1.$
\end{mathitem}

Let $G_1$ be an induced subgraph of $G(k)$ which is isomorphic to $H(k),$ and we name the vertices of $G_1$ accordingly as in $H(k).$
Denote $S=V(G_1)\setminus \{v_1, u_1\}$ and let $G_2=G(k)-S.$ By symmetry, to show that $G(k)$ is pancyclic, it suffices to show that every edge in $\Omega$ is
pancyclic in $G(k).$

It follows from (i) that
\vspace{-5mm}
\begin{itemize}
\item[({\it vii})] for every integer $p\in [3(k-1),\, (6k-5)(k-1)],$ $G_2$ contains a $(v_1,u_1)$-path of length $p.$
\end{itemize}
\vspace{-5mm}
By (ii) and (vii) we deduce that for every integer $p\in [9k-8,\, 6k^2-5k],$ every edge in $\Omega$ lies in a $p$-cycle.

By (iii) and (vii), for every index $i\in [1,\, 3k-4]$ the edge $v_iv_{i+1}$ lies in a $p$-cycle for every integer $p$ in
$$
[3,\, 6k-4]\cup [3k-i+1+3(k-1),\, 3k-i+1+(6k-5)(k-1)]=[3,\, 6k-4]\cup [6k-i-2,\, 6k^2-8k-i+6].
$$
Since $k\ge 3$ and $i\in [1,\, 3k-4],$ we have  $6k-i-2\le 6k-3$ and $6k^2-8k-i+6>9k-8.$ Hence every edge $v_iv_{i+1}$
is pancyclic in $G(k).$

By (iv) and (vii), the edge $vu_{3k-3}$ lies in a $p$-cycle for every integer $p$ in
$$
[3,\, 6k-4]\cup [4+3(k-1),\, 4+(6k-5)(k-1)]=[3,\, 6k-4]\cup [3k+1,\, 6k^2-11k+9].
$$
Since $3k+1\le 6k-3$ and $6k^2-11k+9>9k-8,$ $vu_{3k-3}$ is pancyclic in $G(k).$

By (v) and (vii), the edge $vu$ lies in a $p$-cycle for every integer $p$ in
$$
[3,\, 3k-1]\cup [3+3(k-1),\, 3+(6k-5)(k-1)]=[3,\, 3k-1]\cup [3k,\, 6k^2-11k+8].
$$
Since $6k^2-11k+8>9k-8,$ $vu$ is pancyclic in $G(k).$

By (vi) and (vii), for every index $i\in [1,\, 3k-3]$ the edge $vv_i$ lies in a $p$-cycle for every integer $p$ in
$$
[3,\, 6k-i-3]\cup [3k-i+1+3(k-1),\, 3k-i+1+(k-1)(6k-5)]=[3,\, 6k-i-3]\cup [6k-i-2, 6k^2-8k-i+6].
$$
Since $6k^2-8k-i+6>9k-8,$ $vv_i$ is pancyclic in $G(k).$ This completes the proof. \hfill$\Box$

\section{The maximum diameter of an edge-pancyclic graph of a given order}

{\bf Theorem 6.} {\it The maximum diameter of an edge-pancyclic graph of order $n$ is $\lfloor 2n/5\rfloor.$  }

{\bf Proof.} Let $G$ be an edge-pancyclic graph of order $n$ and diameter $d.$ Since $G$ is edge-pancyclic, it is $2$-connected. Hence
$\delta(G)\ge 2.$ We claim that if $n\ge 4,$ then $\delta(G)\ge 3.$ To the contrary, suppose that $G$ has a vertex $v$ of degree $2$ and let $N(v)=\{v_1, v_2\}.$
If $v_1$ and $v_2$ are adjacent, then the edge $v_1v_2$ is not contained in a Hamilton cycle, a contradiction. If $v_1$ and $v_2$ are nonadjacent, then the edge $vv_1$ is not contained in a triangle, a contradiction again.

If $3\le n\le 5,$ the inequality $d\le \lfloor 2n/5\rfloor$ can be checked directly. It is easy to see that if $H$ is any hamiltonian graph of order $n\ge 6$ with $\delta(H)\ge 3,$ then the diameter of $H$ is at most $\lfloor n/2\rfloor -1.$  Thus if $6\le n\le 9,$  we have $d\le \lfloor n/2\rfloor -1=\lfloor 2n/5\rfloor. $

Next we assume that $n\ge 10.$ Then $\lfloor 2n/5\rfloor\ge 4.$ Thus if $d\le 4,$ the inequality $d\le \lfloor 2n/5\rfloor$ holds trivially. We further assume that
$d\ge 5.$

Choose a peripheral vertex $x$ of $G;$ i.e., the eccentricity of $x$ is $d.$ Consider the distance layers of $x$ in $G:$
$$
V_i\triangleq \{v |\,d(x,v)=i,\,\,\,v\in V(G)\},\quad i=0,1,2,\ldots,d.
$$
Since $G$ is edge-pancyclic, it is $2$-connected. Hence $|V_i|\ge 2$ for $1\le i\le d-1.$

We assert that
$$
|V_1|\ge 3, \eqno (4)
$$
\vspace{-9mm}
$$
|V_{d-1}|+|V_d|\ge 4, \eqno(5)
$$
\vspace{-9mm}
$$
|V_i|+|V_{i+1}|\ge 5, \quad i=1,2,\ldots,d-2.\eqno(6)
$$

The inequality (4) follows from $\delta(G)\ge 3$ proved above.

To prove (5), to the contrary, suppose that $|V_{d-1}|+|V_d|\le 3.$ Then $|V_{d-1}|=2$ and  $|V_d|=1.$ However, the fact that $\delta(G)\ge 3$ implies
$|V_{d-1}|\ge 3,$ a contradiction.

To prove (6), to the contrary, suppose that $|V_i|+|V_{i+1}|\le 4.$ Then $|V_i|$=$|V_{i+1}|=2.$ Let $e$ be an edge with one endpoint in $V_{i}$ and the other in $V_{i+1}.$
Since $e$ lies in a triangle $T,$ one edge $f=yz$ of $T$ lies in $V_{i}$ or in $V_{i+1}.$ Hence one of $V_{i}$ and $V_{i+1}$ is equal to  $\{y, z\}.$ It follows that
$\{y, z\}$ is a vertex cut of $G,$ implying that $f$ does not lie in a Hamilton cycle. This contradicts the condition that $G$ is edge-pancyclic.

If $d$ is odd, we have
\begin{align*}
n=\sum_{i=0}^d |V_i|&=|V_0|+|V_1|+(|V_{d-1}|+|V_d|)+\sum_{j=1}^{(d-3)/2}(|V_{2j}|+|V_{2j+1}|)\\
                    &\ge 1+3+4+((d-3)/2)\cdot 5\\
                    &=(5d+1)/2,
\end{align*}
yielding $d\le (2n-1)/5.$

If $d$ is even, we have
\begin{align*}
n=\sum_{i=0}^d |V_i|&=|V_0|+|V_1|+|V_{d-2}|+(|V_{d-1}|+|V_d|)+\sum_{j=1}^{(d-4)/2}(|V_{2j}|+|V_{2j+1}|)\\
                    &\ge 1+3+2+4+((d-4)/2)\cdot 5\\
                    &=5d/2,
\end{align*}
yielding $d\le 2n/5.$ Since $d$ is an integer, $d\le \lfloor 2n/5\rfloor.$

Now we construct extremal graphs to show that for every integer $n\ge 3,$ the upper bound $\lfloor 2n/5\rfloor$ can be attained, and hence it is in fact the maximum diameter.
The triangle is the extremal graph for the order $n=3.$ For $4\le n\le 7,$ the wheel graph $W_n$ is an extremal graph. Two extremal graphs for the orders $n=8, 9$ are depicted in
Figure 4.
\begin{figure}[h]
\centering
\includegraphics[width=0.7\textwidth]{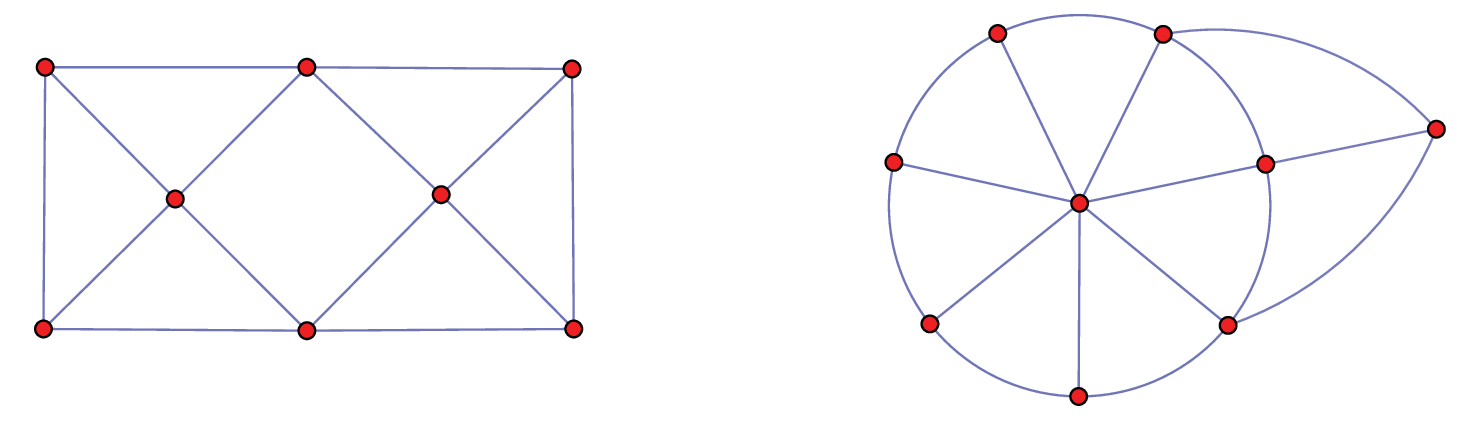}
\caption{Extremal graphs of orders $8$ and $9$}
\end{figure}

Next assume $n\ge 10.$ For pairwise vertex-disjoint graphs $G_1, G_2, \ldots, G_s,$ we denote by  $G_1\vee G_2\vee\cdots\vee G_s$ the {\it sequential join} of $G_1, G_2, \ldots, G_s,$
which is obtained from the disjoint union $G_1+G_2+\cdots+G_s$ by adding edges joining every vertex in $G_i$ to every vertex in $G_{i+1}$ for $i=1,2,\ldots,s-1.$

We distinguish five cases according to the remainder of $n$ divided by $5,$ and construct an extremal graph $Q_n$ for every order $n.$ The building blocks are the triangle
$C_3$ and the empty graph $\overline{K_2}$ of order $2.$

If $n=5k,$ we have $\lfloor 2n/5\rfloor=2k.$ Let $G_1=G_{2k+1}=K_1,$ $G_{2j}=C_3$ for $j=1,\ldots,k$ and $G_{2j+1}=\overline{K_2}$ for $j=1,\ldots,k-1.$
Define $Q_n=G_1\vee G_2\vee\cdots\vee G_{2k+1}.$

If $n=5k+1,$ we have $\lfloor 2n/5\rfloor=2k.$ Let $G_1=K_1,$ $G_{2k+1}=K_2,$ $G_{2j}=C_3$ for $j=1,\ldots,k$ and $G_{2j+1}=\overline{K_2}$ for $j=1,\ldots,k-1.$
Define $Q_n=G_1\vee G_2\vee\cdots\vee G_{2k+1}.$

If $n=5k+2,$ we have $\lfloor 2n/5\rfloor=2k.$ Let $G_1=G_{2k+1}=K_2,$ $G_{2j}=C_3$ for $j=1,\ldots,k$ and $G_{2j+1}=\overline{K_2}$ for $j=1,\ldots,k-1.$
Define $Q_n=G_1\vee G_2\vee\cdots\vee G_{2k+1}.$

If $n=5k+3,$ we have $\lfloor 2n/5\rfloor=2k+1.$ Let $G_1=G_{2k+2}=K_1,$ $G_{2k+1}=G_{2j}=C_3$ for $j=1,\ldots,k$ and $G_{2j+1}=\overline{K_2}$ for $j=1,\ldots,k-1.$
Define $Q_n=G_1\vee G_2\vee\cdots\vee G_{2k+2}.$

If $n=5k+4,$ we have $\lfloor 2n/5\rfloor=2k+1.$ Let $G_1=K_1,$ $G_{2k+2}=K_2,$ $G_{2k+1}=G_{2j}=C_3$ for $j=1,\ldots,k$ and $G_{2j+1}=\overline{K_2}$ for $j=1,\ldots,k-1.$
Define $Q_n=G_1\vee G_2\vee\cdots\vee G_{2k+2}.$

It is not hard to verify that $Q_n$ is an edge-pancyclic graph of order $n$ and diameter $\lfloor 2n/5\rfloor.$ We omit the tedious details. This completes the proof.
\hfill$\Box$

Finally we pose the following problem.

{\bf Problem.} Determine the minimum size of an edge-pancyclic graph of order $n.$

\vskip 5mm
{\bf Acknowledgement.} This research  was supported by the NSFC grant 12271170 and Science and Technology Commission of Shanghai Municipality
 grant 22DZ2229014.

\end{document}